\documentclass[12pt]{article}
\usepackage{latexsym}
\usepackage{amssymb}
\usepackage{graphicx}
\usepackage{cite}
\usepackage{color}
\usepackage{ifpdf}

\newtheorem{Theorem}{Theorem}[part]
\newtheorem{Definition}{Definition}[part]
\newtheorem{Proposition}{Proposition}[part]

\newtheorem{Corollary}{Corollary}[part]

\newtheorem{Conjecture}{Conjecture}

\def \ep{\hbox{ }\hfill$\Box$}

\begin{document}
\title{Regular Uniform Hypergraphs, $s$-Cycles, $s$-Paths and Their largest Laplacian H-Eigenvalues\thanks{This
work was supported by the Hong Kong Research Grant
Council (Grant No. PolyU 501909, 502510, 502111 and 501212) and NSF of China (Grant No. 11231004).}}

\author{
Liqun Qi \thanks{Email: maqilq@polyu.edu.hk. Department of Applied
Mathematics, The Hong Kong Polytechnic University, Hung Hom,
Kowloon, Hong Kong.}, \hspace{4mm} Jia-Yu Shao \thanks{Email:
jyshao@tongji.edu.cn. Department of Mathematics, Tongji University,
Shanghai, China.}\hspace{4mm} Qun Wang\thanks{Email:
wangqun876@gmail.com. Department of Applied Mathematics, The Hong
Kong Polytechnic University, Hung Hom, Kowloon, Hong Kong.}}

\date{\today}
\maketitle

%---------------------------------------------------------------------------------Abstract
\begin{abstract}
In this paper, we show that the largest signless Laplacian
H-eigenvalue of a connected $k$-uniform hypergraph $G$, where $k \ge
3$, reaches its upper bound $2\Delta(G)$, where $\Delta(G)$ is the
largest degree of $G$, if and only if $G$ is regular.   Thus the
largest Laplacian H-eigenvalue of $G$, reaches the same upper bound,
if and only if $G$ is regular and odd-bipartite.   We show that an
$s$-cycle $G$, as a $k$-uniform hypergraph, where $1 \le s \le k-1$,
is regular if and only if there is a positive integer $q$ such that
$k=q(k-s)$. We show that an even-uniform $s$-path and an
even-uniform non-regular $s$-cycle are always odd-bipartite.  We
prove that a regular $s$-cycle $G$ with $k=q(k-s)$ is odd-bipartite
if and only if $m$ is a multiple of $2^{t_0}$, where $m$ is the
number of edges in $G$, and $q = 2^{t_0}(2l_0+1)$ for some integers
$t_0$ and $l_0$. We identify the value of the largest signless
Laplacian H-eigenvalue of an $s$-cycle $G$ in all possible cases.
When $G$ is odd-bipartite, this is also its largest Laplacian
H-eigenvalue.  We introduce supervertices for hypergraphs, and show
the components of a Laplacian H-eigenvector of an odd-uniform
hypergraph are equal if such components corresponds vertices in the
same supervertex, and the corresponding Laplacian H-eigenvalue is
not equal to the degree of the supervertex.   Using this property,
we show that the largest Laplacian H-eigenvalue of an odd-uniform
generalized loose $s$-cycle $G$ is equal to $\Delta(G)=2$. We also
show that the largest Laplacian H-eigenvalue of a $k$-uniform tight
$s$-cycle $G$ is not less than $\Delta(G)+1$, if the number of
vertices is even and $k=4l+3$ for some nonnegative integer $l$.

\vspace{2mm} \noindent {\bf Key words:}\hspace{2mm} Regular uniform
hypergraph, loose cycle, loose path, tight cycle, tight path,
H-eigenvalue, Laplacian.  \vspace{1mm}

\noindent {\bf MSC (2010):}\hspace{2mm}
05C65; 15A18
\end{abstract}

%%%%%%%%%%%%%%%%%%%%%%%%%%%%%%%%%%%%%%%%%%%%%%%%%%%%%%%%%%%%%%%%%%%%%%%%%%%%%%%%%%%%%%%%%%%
\section{Introduction}
\setcounter{Theorem}{0} \setcounter{Proposition}{0}
\setcounter{Corollary}{0} \setcounter{Lemma}{0}
\setcounter{Definition}{0} \setcounter{Remark}{0}
\setcounter{Conjecture}{0}  \setcounter{Example}{0} \hspace{4mm}

Let $k \ge 2$ and $n \ge k$.  A $k$-uniform hypergraph $G = (V, E)$
has vertex set $V$, which is labeled as $[n]=\{1,\ldots,n\}$, and
edge set $E$. By $k$-uniformity, we mean that for every edge $e\in
E$, the cardinality $|e|$ of $e$ is equal to $k$.  If $k=2$, we have
an ordinary graph.

The largest Laplacian eigenvalue of a graph plays an important role
in spectral graph theory \cite{bh11, z07}.   A natural definition
for the Laplacian and signless Laplacian tensors of a $k$-uniform
hypergraph $G$, where $k \ge 3$, was introduced in \cite{q13}.  It
was shown there that the largest Laplacian H-eigenvalue of $G$ is
always less than or equal to the largest signless Laplacian
H-eigenvalue of $G$, while the latter is always less than or equal
to $2\Delta$, where $\Delta$ is the largest degree of $G$.   In
\cite{hq13}, the odd-bipartite hypergraph was introduced.   In
\cite{hqx13}, it was proved that the largest Laplacian H-eigenvalue
of a connected $k$-uniform hypergraph $G$ is equal to its largest
signless Laplacian H-eigenvalue if and only if $G$ is odd-bipartite.
This result extended the classical result in spectral graph theory
\cite{bh11, z07}.

In this paper, we show that the largest signless Laplacian
H-eigenvalue of a connected $k$-uniform hypergraph $G$, where $k \ge
3$, reaches its upper bound $2\Delta$, if and only if $G$ is
regular.   Thus, the largest Laplacian H-eigenvalue of $G$, reaches
the same upper bound, if and only if $G$ is regular and
odd-bipartite.

We then turn our attention to $s$-paths and $s$-cycles.

Researchers in hypergraph theory have studied loose cycles, loose
paths, tight cycles and tight paths extensively \cite{f10, fkl12,
gss08, hlprrss06, kkmo11, kmo10, ko06, mors13, p13}.

Let $G = (V, E)$ be a $k$-uniform hypergraph.   Suppose $1 \le s \le
k-1$.   According to \cite{p13}, if $V = \{ i : i \in [s+m(k-s)] \}$
such that $\{1+j(k-s), \cdots, s+(j+1)(k-s)\}$ is an edge of $G$ for
$j = 0, \cdots, m-1$, then $G$ is called an $s$-path.   In
\cite{mors13}, $G$ is called a {\sl loose path} if $s=1$, and a {\sl
tight path} if $s = k-1$.   In \cite{p13}, $G$ is also called a
loose path for $2 \le s \le {k \over 2}$ and a tight path for ${k
\over 2} < s \le k-2$.  To avoid confusion, in these two cases, as
in \cite{hqs13}, we call $G$ a {\sl generalized loose path} and a
{\sl generalized tight path} respectively. According to \cite{p13},
if $V = \{ i : i \in [m(k-s)] \}$ such that $\{1+j(k-s), \cdots,
s+(j+1)(k-s)\}$ is an edge of $G$ for $j = 0, \cdots, m-1$, where
vertices $m(k-s)+j \equiv  j$ for any $j$, then $G$ is called an
$s$-cycle.   According to \cite{f10, fkl12, gss08, hlprrss06,
kkmo11, kmo10, ko06, mors13}, if $s=1$, $G$ is called a {\sl loose
cycle}, if $s = k-1$, $G$ is called a {\sl tight cycle}. We call $G$
a {\sl generalized loose cycle} for $2 \le s \le {k \over 2}$, and a
{\sl generalized tight cycle} for ${k \over 2} < s \le k-2$.    For
an $s$-cycle, in this paper, we assume that $n \ge 2k-s$.  In this
way, each pair of consecutive edges in the $s$-cycle will have
exactly $s$ common vertices.  In the next section, we will discuss
this in details.

We show in this paper that an $s$-cycle $G$, as a $k$-uniform
hypergraph, where $1 \le s \le k-1$,  is regular if and only if $k$
is a multiple of $k-s$.

The Laplacian H-eigenvalues of loose paths and loose cycles were
studied in \cite{hqs13, hqx13}.   In \cite{hqs13}, power hypergraphs
and cored hypergraphs were introduced.   Loose paths and loose
cycles are power hypergraphs. Power hypergraphs are cored
hypergraphs.   Even-uniform cored hypergraphs are odd-bipartite.  As
cycles are symmetric, their largest signless Laplacian H-eigenvalues
can be identified directly. Thus, the largest Laplacian
H-eigenvalues of odd-bipartite cycles can be identified directly. In
\cite{hqx13}, the largest Laplacian H-eigenvalue of an even-uniform
loose cycle was identified directly.

According to \cite{q13}, the largest Laplacian H-eigenvalue of
$k$-uniform hypergraph is always greater than or equal to the
largest degree of that $k$-uniform hypergraph.  By \cite{hqx13},
when $k$ is even, equality cannot hold, but when $k$ is odd,
equality may hold in certain cases.   It was proved in \cite{hqs13}
that equality holds for odd-uniform loose paths and loose cycles.

It was observed in \cite{hqs13} that if $2 \le s < {k \over 2}$,
then an $s$-path or an $s$-cycle is a cored hypergraph, but not a
power hypergraph in general.

These results raised several questions.  First, if $k$ is even and
${k \over 2} \le s \le k-1$, are some $s$-paths and $s$-cycles still
odd-bipartite, though they are not cored hypergraphs?  Second, can
we identify the largest Laplacian H-eigenvalues of even-uniform
odd-bipartite $s$-cycles directly?  Third, when $k$ is odd and $2
\le s \le k-1$, are the largest H-eigenvalues of $s$-paths and
$s$-cycles equal to the corresponding largest degrees?   We will
study these questions in this paper.

We give some basic definitions in the next section.

In Section 3, we prove the result about regular uniform hypergraphs
mentioned before, and identify regular $s$-cycles.

In Section 4, we show that if $k$ is even, all the $s$-paths and all
the non-regular $s$-cycles are odd-bipartite.   We prove that a
regular $s$-cycle $G$ with $k=q(k-s)$ is odd-bipartite if and only
if $m$ is a multiple of $2^{t_0}$, where $m$ is the number of edges
in $G$, and $q = 2^{t_0}(2l_0+1)$ for some integers $t_0$ and $l_0$.

In \cite{hqs13, hqx13}, several classes of hypergraphs were shown to
be odd-bipartite.   But only in this paper, some regular $s$-cycles
are shown to be not odd-bipartite.   To show that a hypergraph is
odd-bipartite, one only needs to give an adequate odd-partition for
the vertex set of that hypergraph.  Such an odd-partition is not
unique in general. To show that a hypergraph is not odd-bipartite,
one needs to prove that there is no such an odd-partition for the
vertex set of that hypergraph.   Hence, in general, it is not a
trivial task to show that a hypergraph is not odd-bipartite.

In Section 5, we identify the largest Laplacian H-eigenvalues of
even-uniform $s$-cycles directly when $2 \le s \le k-1$.  These
include all the even-uniform non-regular $s$-cycles, and those
odd-bipartite regular $s$-cycles.

We introduce supervertices for hypergraphs in Section 6, and show
there that the components of an H-eigenvector of an odd-uniform
hypergraph are equal if such components corresponds vertices in the
same supervertex, and the corresponding Laplacian H-eigenvalue is
not equal to the degree of the supervertex.   Using this property,
in Section 7, we show that the largest H-eigenvalue of an
odd-uniform generalized loose $s$-cycle is equal to $2$, the maximum
degree of that $s$-cycle.

In Section 8, we show that the largest Laplacian H-eigenvalue of a
$k$-uniform tight $s$-cycle is at least $k+1$, if the number of
vertices is even and $k=4l+3$ for some nonnegative integer $l$. Note
that in this case $\Delta = k$.    We show that equality holds here
if $k=3, s=2$ and $n=4$. When $k=3, s=2$, and $n \ge 5$, we show
that the largest Laplacian H-eigenvalue is no more than $4.5$.

Some final remarks are made in Section 9.

%%%%%%%%%%%%%%%%%%%%%%%%%%%%%%%%%%%%%%%%%%%%%%%%%%%%%%%%%%
\section{Preliminaries}\label{sec-p}
\setcounter{Theorem}{0} \setcounter{Proposition}{0}
\setcounter{Corollary}{0} \setcounter{Lemma}{0}
\setcounter{Definition}{0} \setcounter{Remark}{0}
\setcounter{Conjecture}{0}  \setcounter{Example}{0} \hspace{4mm}
%%%%%%%%%%%%%%%%%%%%%%%%%%%%%%%%%%%%%%%%%%

Let $\mathbb R$ be the field of real numbers, $\mathbb R^n$ the
$n$-dimensional real space, and $\mathbb R^n_+$ the nonnegative
orthant of $\mathbb R^n$. For integers $k\geq 3$ and $n\geq 2$, a
real tensor $\mathcal T=(t_{i_1\ldots i_k})$ of order $k$ and
dimension $n$ refers to a multidimensional array (also called
hypermatrix) with entries $t_{i_1\ldots i_k}$ such that
$t_{i_1\ldots i_k}\in\mathbb{R}$ for all $i_j\in[n]:=\{1,\ldots,n\}$
and $j\in[k]$. Tensors are always referred to $k$-th order real
tensors in this paper, and the dimensions will be clear from the
content. Given a vector $\mathbf{x}\in \mathbb{R}^{n}$, ${\cal
T}\mathbf{x}^{k-1}$ is defined as an $n$-dimensional vector such
that its $i$-th element is
$\sum\limits_{i_2,\ldots,i_k\in[n]}t_{ii_2\ldots i_k}x_{i_2}\cdots
x_{i_k}$ for all $i\in[n]$. Let ${\cal I}$ be the identity tensor of
appropriate dimension, e.g., $i_{i_1\ldots i_k}=1$ if and only if
$i_1=\cdots=i_k\in [n]$, and zero otherwise when the dimension is
$n$. The following definition was introduced in \cite{q05}.
%---------------------------------------------------------------Definition eigenvalue
\begin{Definition}\label{def-00}
Let $\mathcal T$ be a $k$-th order $n$-dimensional real tensor. For
some $\lambda\in\mathbb{R}$, if polynomial system $\left(\lambda
{\cal I}-{\cal T}\right)\mathbf{x}^{k-1}=0$ has a solution
$\mathbf{x}\in\mathbb{R}^n\setminus\{0\}$, then $\lambda$ is called
an H-eigenvalue and $\mathbf x$ an
H-eigenvector.
\end{Definition}
H-eigenvalues are real numbers, by Definition \ref{def-00}. By
\cite{q05}, we have that the number of H-eigenvalues of a real
tensor is finite. By \cite{q13}, we have that all the tensors
considered in this paper have at least one H-eigenvalue. Hence, we
can denote by $\lambda(\mathcal T)$ as the largest H-eigenvalue of a
real tensor $\mathcal T$.

For a subset $S\subseteq [n]$, we denoted by $|S|$ its cardinality.
%and $\mbox{supp}(\mathbf x):=\{i\in[n]\;|\;x_i\neq 0\}$ is the {\em
%support} of $\mathbf x$.

Consider a $k$-uniform hypergraph $G = (V, E)$ with vertex set $V$,
which is labeled as $[n]=\{1,\ldots,n\}$, and edge set $E$. For a
subset $S\subset [n]$, we denote by $E_S$ the set of edges $\{e\in
E\;|\;S\cap e\neq\emptyset\}$. For a vertex $i\in V$, we simplify
$E_{\{i\}}$ as $E_i$. It is the set of edges containing the vertex
$i$, i.e., $E_i:=\{e\in E\;|\;i\in e\}$. The cardinality $|E_i|$ of
the set $E_i$ is defined as the {\em degree} of the vertex $i$,
which is denoted by $d_i$. If two vertices $i$ and $j$ are in the
same edge $e$, then we denote $i \sim j$.  Two different vertices
$i$ and $j$ are {\em connected} to each other (or the pair $i$ and
$j$ is connected), if there is a sequence of edges
$(e_1,\ldots,e_m)$ such that $i\in e_1$, $j\in e_m$ and $e_r\cap
e_{r+1}\neq\emptyset$ for all $r\in[m-1]$. A hypergraph is called
{\em connected}, if every pair of vertices of $G$ is connected. %Let
%$S\subseteq V$, the hypergraph with vertex set $S$ and edge set
%$\{e\in E\;|\;e\subseteq S\}$ is called the {\em sub-hypergraph} of
%$G$ induced by $S$. We will denote it by $G_S$.
A hypergraph is {\em regular} if $d_1=\cdots=d_n=d$. A hypergraph
$G=(V,E)$ is {\em complete} if $E$ consists of all the possible
edges. In this case, $G$ is regular of degree $d={n-1\choose k-1}$.

 For the sake of simplicity, we mainly consider connected hypergraphs in the subsequent analysis. By the techniques in \cite{q13,hq13}, the conclusions on connected hypergraphs can be easily generalized to general hypergraphs.

The following definition for the Laplacian tensor and signless
Laplacian tensor was proposed by Qi \cite{q13}.
%-----------------------------------------
\begin{Definition}\label{def-l}
Let $G=(V,E)$ be a $k$-uniform hypergraph. The {\em adjacency tensor} of $G$ is defined as the $k$-th order $n$-dimensional tensor $\mathcal A$ whose $(i_1 \ldots i_k)$-entry is:
\begin{eqnarray*}
a_{i_1 \ldots i_k}:=\left\{\begin{array}{cl}\frac{1}{(k-1)!}&if\;\{i_1,\ldots,i_k\}\in E,\\0&\mbox{otherwise}.\end{array}\right.
\end{eqnarray*}
Let $\mathcal D$ be a $k$-th order $n$-dimensional diagonal tensor with its diagonal element $d_{i\ldots i}$ being $d_i$, the degree of vertex $i$, for all $i\in [n]$. Then $\mathcal L:=\mathcal D-\mathcal A$ is the {\em Laplacian tensor} of the hypergraph $G$, and $\mathcal Q:=\mathcal D+\mathcal A$ is the {\em signless Laplacian tensor} of the hypergraph $G$.
\end{Definition}

By \cite{q13}, zero is always the smallest H-eigenvalue of $\mathcal
L$, and we have $\lambda(\mathcal L) \le \lambda(\mathcal Q) \le
2\Delta$, where $\Delta$ is the maximum degree of $G$.   For ${\cal
T} = {\cal L}$, the polynomial system $\left(\lambda {\cal I}-{\cal
T}\right)\mathbf{x}^{k-1}=0$ in Definition \ref{def-00} has the form
\begin{equation} \label{e-1}
\lambda x_i^{k-1} = d_ix_i^{k-1} - \sum_{\{i, i_2, \cdots, i_m\} \in
E} x_{i_2}\cdots x_{i_k}, \ \ {\rm for} \ i \in [n].
\end{equation}
For ${\cal T} = {\cal Q}$, the polynomial system $\left(\lambda
{\cal I}-{\cal T}\right)\mathbf{x}^{k-1}=0$ in Definition
\ref{def-00} has the form
\begin{equation} \label{e-2}
\lambda x_i^{k-1} = d_ix_i^{k-1} + \sum_{\{i, i_2, \cdots, i_m\} \in
E} x_{i_2}\cdots x_{i_k}, \ \ {\rm for} \ i \in [n].
\end{equation}

In the following, we define cored hypergraphs.
%-----------------------------------------
\begin{Definition}\label{def-cor}
Let $G=(V,E)$ be a $k$-uniform hypergraph. If for every edge $e\in
E$, there is a vertex $i_e\in e$ such that the degree of the vertex
$i_e$ is one, then $G$ is called a {\em cored hypergraph}. A vertex
with degree one is called a {\em core vertex}, and a vertex with
degree larger than one is called an {\em intersection vertex}.
\end{Definition}

The notion of odd-bipartite even-uniform hypergraphs was introduced
in \cite{hq13}.
%-----------------------------------------
\begin{Definition}\label{def-bi-odd}
Let $G=(V,E)$ be a $k$-uniform hypergraph. Then $G$ is called {\em
odd-bipartite} if $k$ is even and either it is trivial (i.e.,
$E=\emptyset$) or there is a partition of the vertex set $V$ as
$V=V_1\cup V_2$ such that $V_1,V_2\neq \emptyset$ and every edge in
$E$ intersects $V_1$ with exactly an odd number of vertices.
\end{Definition}

In the introduction, we claim that $n$, the number of vertices in an
$s$-cycle needs to satisfy the condition $n \ge 2k-s$ such that each
pair of consecutive edges has exactly $s$ common vertices.   We now
discuss this in details below.

\begin{Proposition} \label{p-1}
Let $G=(V,E)$ be a k-uniform $s$-cycle with $n$ vertices and $m$
edges, where $n=m(k-s)$ and $1\le s \le k-1$. Then each pair of
consecutive edges of G contains exactly $s$ common vertices if and
only if $n\ge 2k-s$.
\end{Proposition}
\noindent {\bf Proof.} By the definition, we may assume that
$V=[n]$, and may agree that $i=n+i$ when $i$ and $n+i$ are both
viewed as vertices of $G$. Also, we have $E=\{e_0,e_1,\cdots,
e_{m-1}\}$, where
$$e_j=\{j(k-s)+1,\cdots, j(k-s)+k\} \qquad  (j=0,1,\cdots, m-1)$$

Necessity. If each pair of consecutive edges of G contains exactly
$s$ common vertices, then $|e_0\cap e_1|=s$. On the other hand, we
have $|e_0\cup e_1|+|e_0\cap e_1|=|e_0| + |e_1|$. So we have
$$n\ge |e_0\cup e_1|= |e_0| + |e_1|- |e_0\cap e_1| =2k-s$$

Sufficiency. Now suppose that $n\ge 2k-s$. Then it is easy to verify
that
$$e_j\cap e_{j+1}=\{(j+1)(k-s)+1,\cdots, (j+1)(k-s)+s\} \qquad  (j=0,1,\cdots, m-2)$$
Since $n\ge 2k-s$ implying $(m-1)(k-s)\ge k$, we can also verify
that
$$e_{m-1}\cap e_0=\{1,\cdots, s\} $$
From these relations we obtain
$$|e_j\cap e_{j+1}|= s \quad  (j=0,1,\cdots, m-2) \quad \mbox {and} \quad |e_{m-1}\cap e_0|=s.$$
  \ep

%%%%%%%%%%%%%%%%%%%%%%%%%%%%%%%%%%%%%%%%%%%%%%%%%%%%%%%%%%%%%%%%%%%%%%%%%%%%%%%
\section{Regular Uniform Hypergraphs and Regular $s$-Cycles}\label{sec-ch0}
\setcounter{Theorem}{0} \setcounter{Proposition}{0}
\setcounter{Corollary}{0} \setcounter{Lemma}{0}
\setcounter{Definition}{0} \setcounter{Remark}{0}
\setcounter{Conjecture}{0}  \setcounter{Example}{0} \hspace{4mm}

We now establish the following theorem for a connected $k$-uniform
hypergraph $G$.

\begin{Theorem} \label {t0}
Suppose that $G = (V, E)$ is a connected $k$-uniform hypergraph with
$k \ge 2$ and maximum degree $\Delta$.   Then $\lambda({\mathcal Q})
= 2 \Delta$ if and only if $G$ is regular. Furthermore,
$\lambda({\mathcal L}) = 2 \Delta$ if and only if $G$ is regular and
odd-bipartite.
\end{Theorem}

\noindent {\bf Proof.}  First, we assume that $G$ is regular. Then,
by \cite[Theorem 3.4]{q13} and \cite{cqz13} (see also \cite[Lemmas
2.2 and 2.3]{hq13}), if we can find a positive H-eigenvector
$\mathbf x\in\mathbb R^n$ of $\mathcal Q$ corresponding to an
H-eigenvalue $\mu$, then $\mu=\lambda(\mathcal Q)$.   Take $x_i = 1$
for every $i$. By (\ref{e-2}), we have $\mu = \Delta + \Delta =
2\Delta$. Thus, $\lambda({\mathcal Q}) = 2 \Delta$.

On the other hand, assume that $\lambda({\mathcal Q}) = 2 \Delta$.
Let $\mathbf x\in\mathbb R^n$ be a positive H-eigenvector of
$\mathcal Q$ corresponding to the H-eigenvalue $2\Delta$.   Assume
that $x_i = \max \{ x_j : j \in [n] \}$.    Then by (\ref{e-2}), we
have
$$2\Delta x_i^{k-1} = d_ix_i^{k-1} + \sum_{\{ i, i_2, \cdots, i_k \} \in E}  x_{i_2}\cdots x_{i_k},$$
where $d_i$ is the degree of vertex $i$. To make this equality hold,
we must have $d_i = \Delta$ and $x_j = x_i$ as long as $j \sim i$.
Applying the same augment for all such $j$ with $j \sim  i$, we have
$d_j = \Delta$ and $x_l = x_j = x_i$ as long as $l \sim j$.   As $G$
is connected, we see that $d_j = \Delta$ for all $j \in V$.  Thus
$G$ is regular.

The last conclusion of this theorem follows from the above
conclusion and \cite[Theorem 5.1]{hqx13}. \ep

\bigskip

Clearly, an $s$-path cannot be regular.   We now consider regular
$s$-cycles.

\begin{Proposition} \label{p0}
Let $G = (V, E)$ be a $k$-uniform $s$-cycle, with $1 \le s \le k-1$,
$k \ge 3$, $V = \{i : i \in [m(k-s)] \}$, such that $\{ 1+j(k-s),
\cdots, s+(j+1)(k-s)\}$ is an edge of $G$ for $j = 0, \cdots, m-1$,
where vertices $m(k-s)+j \equiv j$ for any $j$. Then $G$ is regular
if and only if $k = q(k-s)$ for some positive integer $q$. In this
case, we have $d_1 = \cdots = d_n = q$, where $n = m(k-s) = | V |$.
\end{Proposition}

\noindent {\bf Proof.}  If $k = q(k-s)$, then we see that $d_1 =
\cdots = d_n = q$.   Hence $G$ is regular in this case. On the other
case, suppose $k = q(k-s) + r$, where $1 \le r < k-s$.
 Then we see that $d_1 = q+1$ and $d_{k-s} = q$.   Thus, $G$ cannot be regular in this case.

 The conclusions of this proposition follow now.  \ep

 Since $1 \le s \le k-1$, we see that $2 \le q \le k$.  For a tight cycle, $s = k-1$, we see
 that $G$ is also regular with $q = k$ in this case.  Thus, we have
 the following corollary.

 \begin{Corollary}
 A tight cycle is regular.
 \end{Corollary}

%%%%%%%%%%%%%%%%%%%%%%%%%%%%%%%%%%%%%%%%%%%%%%%%%%%%%%%%%%%%%%%%%%%%%%%%%%%%%%%
\section{Odd-Bipartite $s$-Paths and $s$-Cycles}\label{sec-ch}
\setcounter{Theorem}{0} \setcounter{Proposition}{0}
\setcounter{Corollary}{0} \setcounter{Lemma}{0}
\setcounter{Definition}{0} \setcounter{Remark}{0}
\setcounter{Conjecture}{0}  \setcounter{Example}{0} \hspace{4mm}

We assume that $k$ is even in this section, as odd-bipartite
hypergraphs are only for even $k$.

%%%%%%%%%%%%%%%%%%%%%%%%%%%%%%%%%%%%%%%%%%%%%%%%%%%%%%
\subsection{Odd-bipartite $s$-Paths}\label{sec-ob-1}

Our first proposition in this section shows that when $k$ is even,
all the $s$-paths are odd-bipartite.

%----------------------------------------------------------
\begin{Proposition}\label{p1}
Assume that $k \ge 4$ is even.   Let $G=(V,E)$ be a $k$-uniform
$s$-path, where $1 \le s \le k-1$.   Then $G$ is odd-bipartite.
\end{Proposition}

\noindent {\bf Proof.}  According the discussion at the beginning of
this paper, we may assume that $V = \{i : i\in [s+m(k-s)] \}$ such
that $\{1+j(k-s), \cdots, s+(j+1)(k-s)\}$ is an edge of $G$ for $j =
0, \cdots, m-1$.    Let $V_1 = \{ k, 2k, \cdots, \}$ and $V_2 = V
\setminus V_1$.   Then we see that $G$ is odd-bipartite as each edge
has exactly one vertex in $V_1$. \ep

%%%%%%%%%%%%%%%%%%%%%%%%%%%%%%%%%%%%%%%%%%%%%%%%%%%%%%
\subsection{Odd-bipartite Non-Regular $s$-Cycles}\label{sec-ob-2}

Our second proposition in this section shows that when $k$ is even,
all the non-regular $s$-cycles are odd-bipartite.

\begin{Proposition}\label{p2}
Assume that $k \ge 4$ is even.   Let $G=(V,E)$ be a $k$-uniform
non-regular $s$-cycle, where $1 \le s \le k-1$.   Then $G$ is
odd-bipartite.
\end{Proposition}

\noindent {\bf Proof.}  When $s=1$, $G$ is a loose cycle, thus a
power hypergraph \cite{hqs13}. When $1 < s < {k \over 2}$, $G$ has
at least one core vertex, thus is a cored hypergraph \cite{hqs13}.
In both cases, $G$ is odd-bipartite as long as $k$ is even, as
observed in \cite{hqs13}.

Now, by Proposition \ref{p0}, the remaining case, after excluding
regular $s$-cycles, are that ${k \over 2} < s < k-2$, $k =
q(k-s)+r$, where $1 \le r \le k-s-1$.   We may assume that $V = \{ i
: i \in [m(k-s)] \}$, and note that vertices $j+m(k-s) \equiv j$ for
all $j$. If $q$ is odd, let $V_1 = \{ i(k-s) : i \in [m] \}$ and
$V_2 = V \setminus V_1$. Then each edge has exactly $q$ vertices in
$V_1$. If $q$ is even, let $V_1 = \{ 1+ (i-1)(k-s) : i \in [m] \}$
and $V_2 = V \setminus V_1$. Then each edge has exactly $q+1$
vertices in $V_1$.  In both cases,  each edge has odd number of
vertices in $V_1$.  Thus, $G$ is odd-bipartite as long as $k$ is
even. \ep

%%%%%%%%%%%%%%%%%%%%%%%%%%%%%%%%%%%%%%%%%%%%%%%%%%%%%%
\subsection{Odd-bipartite Regular $s$-Cycles}\label{sec-ob-3}

We now give a sufficient and necessary condition for a regular
$s$-cycle to be odd-bipartite.

%----------------------------------------------------------
\begin{Theorem}\label{t1}
Let $G=(V,E)$ be a k-uniform $s$-cycle with $n$ vertices and $m$
edges, where $n=m(k-s)$, $k$ is even and $1\le s \le k-1$. Assume
that there is an integer $q$ such that $k=q(k-s)$ (Thus $G$ is
regular by Proposition \ref{p0}). Write $q=2^{t_0}(2l_0+1)$ for some
nonnegative integers $t_0$ and $l_0$. Then $G$ is odd-bipartite if
and only if $m$ is a multiple of $2^{t_0}$.
\end{Theorem}

\noindent {\bf Proof.}  We assume that $V=Z_n$ (the set of integers
modulo n). Namely, we agree that $i=n+i$, when $i$ and $n+i$ are
both viewed as vertices of $G$.

Also, the edges $e_0,e_1,\cdots,e_{m-1}$ of $G$ are as follows:
\begin{equation} \label{4.1}
e_j=\{j(k-s)+1,\cdots , j(k-s)+k\}, \qquad  (j=0,1,\cdots, m-1),
\end{equation}
where each edge consists of k cyclicly consecutive vertices of $G$.

\vskip 0.18cm

Sufficiency. Suppose that $m=2^{t_0}p_0$. Let $q_0=2^{t_0}(k-s)$.
Then we have $n=m(k-s)=p_0q_0$. Let
$$V_1=\{q_0, 2q_0\cdots , p_0q_0\} $$
be the set of all the multiples of $q_0$ in the set $Z_n$.

Since $k=q_0(2l_0+1)$ and $n=p_0q_0$ are both multiples of $q_0$, we
see that each set of  k cyclicly consecutive elements in $Z_n$
contains exactly $\frac {k} {q_0}=2l_0+1$ elements which are
multiples of $q_0$. Thus each edge of $G$ contains exactly $2l_0+1$
vertices in $V_1$.  Hence $G$ is odd-bipartite.

\vskip 0.18cm

Necessity. We write $a\sim b$ if the two integers $a$ and $b$ have
the same parity. Let
\begin{equation} \label{4.2}
X_j=\{j(k-s)+1,\cdots , (j+1)(k-s)\}, \qquad  (j=0,1,\cdots, m-1).
\end{equation}
Then we have  $|X_j|=k-s$ and $V=Z_n=\cup_{j=0}^{m-1}X_j$. Also we
agree that $X_{m+j}=X_j$ (as subsets of $Z_n$).

\vskip 0.18cm

By comparing (\ref{4.1}) and (\ref{4.2}) we have
\begin{equation} \label{4.3}
e_i=\bigcup_{j=i}^{q-1+i}X_j,  \qquad  (i=0,1,\cdots, m-1).
\end{equation}

Now suppose that $G$ is odd-bipartite with the bipartition
$(V_1,V_2)$. Let
$$a_j=|X_j\cap V_1| \qquad  \mbox {and} \qquad b_j=|e_j\cap V_1|,  \qquad  (j\equiv 0,1,\cdots, m-1 \  \ mod \ m).$$
Then by the definition of odd-bipartition, all
$b_0,b_1,\cdots,b_{m-1}$ are odd.

\vskip 0.18cm

By (\ref{4.3}) we also have
$$e_i\cap V_1=\left (\bigcup_{j=i}^{q-1+i}X_j \right )\cap V_1=\bigcup_{j=i}^{q-1+i}(X_j\cap V_1)  $$
which implies that
$$b_i=\sum_{j=i}^{q-1+i}a_j  \qquad  (i\equiv 0,1,\cdots, m-1 \  \ mod \ m) $$
From this we have
\begin{equation} \label{4.5}
a_{q+i}-a_i=\sum_{j=i+1}^{q+i}a_j - \sum_{j=i}^{q-1+i}a_j
=b_{i+1}-b_i  \sim 0,  \qquad (i\equiv 0,1,\cdots, m-1 \  \ mod \
m).
\end{equation}

On the other hand, since $X_{m+i}=X_i$ we also have
\begin{equation} \label{4.6}
a_{m+i}= a_i,  \qquad (i\equiv 0,1,\cdots, m-1 \  \ mod \ m).
\end{equation}

Let $g=gcd (m,q)$ be the greatest common divisor of $m$ and $q$.
Then $g=cm+dq$ for some integers $c$ and $d$. So by (\ref{4.5}) and
(\ref{4.6}) we have
\begin{equation} \label{4.7}
a_{g+i}\sim a_i, \qquad (i\equiv 0,1,\cdots, m-1 \  \ mod \ m).
\end{equation}

Now let $q'= q/g$. Then $q=gq'$, so by (\ref{4.7}) we have
$$b_0=\sum_{j=0}^{q-1}a_j =\sum_{t=0}^{q'-1}\sum_{j=0}^{g-1}a_{tg+j}\sim \sum_{t=0}^{q'-1}\sum_{j=0}^{g-1}a_{j}
=q'\sum_{j=0}^{g-1}a_{j}$$

Since $b_0$ is odd, $q'$ is also odd, which implies that $g$ is a
multiple of $2^{t_0}$. Thus $m$ is also a multiple of $2^{t_0}$. \ep

Figure 1 indicates an odd-bipartite regular $3$-cycle with $k =
2(k-s) = 6$ and $m = 4$.   We see that $G$ is odd-bipartite with
$V_1 = \{ 6, 12\}$ and $V_2 = V \setminus V_1$.

%%%%%%%%%%%%%%%%%%%%%%%%%%%%%%%%%%%
\begin{figure}[htbp]
\centering
\includegraphics[width=1.6in]{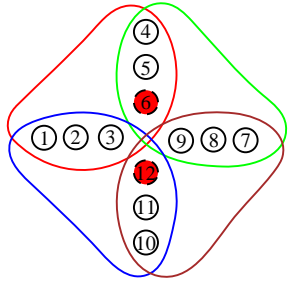}
\caption{An odd-bipartite regular $3$-cycle with $k = 2(k-s) = 6$
and $m = 4$.}
\end{figure}

%Figure 2 indicates an odd-bipartite generalized loose $2$-cycle with
%$k=6$ and $m=2$.   We see that $G$ is odd-bipartite with $V_1 = \{
%3, 8\}$ and $V_2 = V \setminus V_1$.

%%%%%%%%%%%%%%%%%%%%%%%%%%%%%%%%%%%
%\begin{figure}[htbp]
%\centering
%\includegraphics[width=1.6in]{figure2.eps}
%\caption{An odd-bipartite generalized loose $2$-cycle with $k=6$ and
%$m=2$.}
%\end{figure}

The smallest non-odd-bipartite even-uniform regular $s$-cycle may be
as follows: $k=4, s=2$ and $m=3$ (thus $n=6$.   Using Matlab, we
find that the Laplacian H-eigenvalues of this $2$-cycle are $0, 1,
2$ and $3$ only.  Thus, $\lambda({\cal L})= 3 < 2\Delta = 4$ in this
example.    This confirms Theorem \ref{t1}.

\begin{Corollary} \label{c1}
In Theorem \ref{t1}, if $k$ is even and $q$ is odd, then $G$ is
odd-bipartite.
\end{Corollary}

If $s = k-1$, then $G$ is a tight cycle.   We have the following
corollary.

\begin{Corollary} \label{c2}
Let $G=(V,E)$ be a $k$-uniform tight cycle, i.e., $s = k-1$. Then
$G$ is regular.  Assume that $k$ is even.   We may write $k =
2^{t_0}(2l_0+1)$ for two nonnegative integers $t_0$ and $l_0$. Then
$G$ is odd-bipartite if and only if $m=n$ is a multiple of
$2^{t_0}$.
\end{Corollary}

%%%%%%%%%%%%%%%%%%%%%%%%%%%%%%%%%%%%%%%%%%%%%%%%%%%%%%%%%%%%%%%%%%%%%%%%%%%%%%%
\section{The Largest Laplacian H-eigenvalue of Odd-Bipartite $s$-Cycles}\label{sec-ph}
\setcounter{Theorem}{0} \setcounter{Proposition}{0}
\setcounter{Corollary}{0} \setcounter{Lemma}{0}
\setcounter{Definition}{0} \setcounter{Remark}{0}
\setcounter{Conjecture}{0}  \setcounter{Example}{0} \hspace{4mm}

In this section, we identify the largest signless Laplacian
H-eigenvalues of $s$-cycles in all possible cases.  When these
$s$-cycles are odd-bipartite, these values are also their largest
Laplacian H-eigenvalues.

%%%%%%%%%%%%%%%%%%%%%%%%%%%%%%%%%%%%%%%%%%%%%%%%%%%%%%
\subsection{$s$-Cycles with Core Vertices}\label{sec-la-1}

This is the case that $1 \le s < {k \over 2}$.  Suppose $G = (V, E)$
is such an $s$-cycle.   Then for each edge, there are $k-2s$ core
vertices, and $2s$ intersection vertices.

By \cite[Theorem 3.4]{q13} and \cite{cqz13} (see also \cite[Lemmas
2.2 and 2.3]{hq13}), if we can find a positive H-eigenvector
$\mathbf x\in\mathbb R^n$ of $\mathcal Q$ corresponding to an
H-eigenvalue $\mu$, then $\mu=\lambda(\mathcal Q)$.

Now we take ${\bf x} \in \Re^n$ be a positive vector with $x_i =
\alpha
> 0$ if $i$ is a core vertex, and $x_j = 1$ if $j$ is an
intersection vertex.  Suppose that $\mathbf x$ is an H-eigenvector
of $\mathcal Q$ corresponding to the H-eigenvalue
$\mu=\lambda(\mathcal Q)$. Note that the degree of a core vertex is
$1$ and the degree of an intersection vertex is $2$.   By
(\ref{e-2}), we would have
\begin{eqnarray*}
\mu = 2 + 2 \alpha^{k-2s}\;\mbox{and}\;\mu \alpha^{k-1}=\alpha^{k-1}
+ \alpha^{k-2s-1},
\end{eqnarray*}
i.e.,
\begin{eqnarray*}
\mu = 2 + 2 \alpha^{k-2s}\;\mbox{and}\;(\mu-1) \alpha^{2s}=1.
\end{eqnarray*}
Eliminating $\mu$, we have
\begin{equation} \label{e1}
f(\alpha) \equiv 2\alpha^k + \alpha^{2s} -1 = 0.
\end{equation}
Since $f(0) = -1 < 0$ and $f(1) = 2 > 0$, (\ref{e1}) has a root
$\alpha_* \in (0, 1)$.   Let $\mu_* = 2 + 2 \alpha_*^{k-2s}$.  Then
$\lambda(\mathcal Q) = \mu_*$.  By \cite{hqs13}, since $G$ is a
cored hypergraph, we have $\lambda(\mathcal L) = \mu_*$ if $k$ is
even.

We conclude this discussion as the following theorem.

\begin{Theorem} \label{t3}
Suppose that $G= (V, E)$ is an $s$-cycle with $k \ge 3$ and $1 \le s
< {k \over 2}$.  Then $\lambda(\mathcal Q) = 2 + 2 \alpha_*^{k-2s}$,
where $\alpha_*$ is the unique root of (\ref{e1}) in $(0, 1)$.  When
$k$ is even, we have $\lambda(\mathcal L) = 2 + 2 \alpha_*^{k-2s}$
too.
\end{Theorem}

Note that in this case $\Delta = 2$ and we have $\Delta = 2 <
\lambda(\mathcal Q) < 2\Delta = 4$.  This confirms Theorem \ref{t0}
and \cite[Corollary 6.2]{q13}.

%%%%%%%%%%%%%%%%%%%%%%%%%%%%%%%%%%%%%%%%%%%%%%%%%%%%%%
\subsection{Regular $s$-Cycles}\label{sec-la-2}

By Proposition \ref{p0} and Theorem \ref{t1}, we have the following
proposition.

\begin{Proposition} \label{p4}
Suppose that $G= (V, E)$ is an $s$-cycle with $k = q(k-s) \ge 3$.
Then $G$ is a regular hypergraph and $\lambda(\mathcal Q) = 2q$.
Assume further that $k$ is even.  If either $q$ is odd or $q =
2^{t_0}(2l_0+1)$ for a positive integer $t_0$ and a nonnegative
integer $l_0$, and $m$ is a multiple of $2^{t_0}$, then we have
$\lambda(\mathcal L) = 2q$ too.
\end{Proposition}

It will be a further research topic to find the value of
$\lambda(\mathcal L)$ for a non-odd-bipartite regular $s$-cycle.

\subsection{Non-Regular Generalized Tight $s$-Cycles}\label{sec-la-3}

In this subsection, we consider a generalized tight $s$-cycle $G=(V,
E)$, which is not a regular $s$-cycle.    We may assume that ${k
\over 2} < s < k-1$ and $k= q(k-s) + r$, where $1 \le r < k-s$. Then
for each edge, there are $(q+1)r$ vertices with degree $q+1$, and
$k-(q+1)r$ vertices with degree $q$.

Again, if we can find a positive H-eigenvector $\mathbf x\in\mathbb
R^n$ of $\mathcal Q$ corresponding to an H-eigenvalue $\mu$, then
$\mu=\lambda(\mathcal Q)$.

Now we take ${\bf x} \in \Re^n$ be a positive vector with $x_i =
\alpha > 0$ if $i$ is a vertex with degree $q$, and $x_j = 1$ if $j$
is a vertex with degree $q+1$.  Suppose that $\mathbf x$ is an
H-eigenvector of $\mathcal Q$ corresponding to the H-eigenvalue
$\mu=\lambda(\mathcal Q)$.  By (\ref{e-2}), we should have
\begin{eqnarray*}
\mu =  q+1 + (q+1) \alpha^{k-(q+1)r}\;\mbox{and}\;\mu
\alpha^{k-1}=q\alpha^{k-1} + q\alpha^{k-(q+1)r-1},
\end{eqnarray*}
i.e.,
\begin{eqnarray*}
\mu = q+1 + (q+1) \alpha^{k-(q+1)r}\;\mbox{and}\;\mu = q +
q\alpha^{-(q+1)r}.
\end{eqnarray*}
Eliminating $\mu$, we have
\begin{equation} \label{e2}
f(\alpha) \equiv (q+1)\alpha^k + \alpha^{(q+1)r} -q = 0.
\end{equation}
Since $f(0) = -q < 0$ and $f(1) = 2 > 0$, (\ref{e2}) has a root
$\alpha_* \in (0, 1)$.   Let $\mu_* = q+1 + (q+1)
\alpha_*^{k-(q+1)r}$. Then $\lambda(\mathcal Q) = \mu_*$.  If $k$ is
even, then $G$ is odd-bipartite.  By \cite{hqx13}, we have
$\lambda(\mathcal L) = \mu_*$ in this case.

We conclude this discussion as the following theorem.

\begin{Theorem} \label{t4}
Suppose that $G= (V, E)$ is an $s$-cycle with $k \ge 3$ and ${k
\over 2} < s < k-1$, $k=q(k-s)+r$, with $1 \le r < k-s$.  Then
$\lambda(\mathcal Q) = q+1 + (q+1) \alpha_*^{k-(q+1)r}$, where
$\alpha_*$ is the unique root of (\ref{e2}) in $(0, 1)$.   When $k$
is even, we have $\lambda(\mathcal L) = q+1 + (q+1)
\alpha_*^{k-(q+1)r}$ too.
\end{Theorem}

In this case $\Delta = q+1$ and we have $\Delta = q+1 <
\lambda(\mathcal Q) < 2\Delta = 2(q+1)$.  This also confirms Theorem
\ref{t0} and \cite[Corollary 6.2]{q13}.

Note that all $s$-cycles are covered by the discussion in these
three subsections.

%%%%%%%%%%%%%%%%%%%%%%%%%%%%%%%%%%%%%%%%%%%%%%%%%%%%%%%%%%%%%%%%%%%%%%%%%%%%%%%
\section{Supervertices}\label{sec-ph1}
\setcounter{Theorem}{0} \setcounter{Proposition}{0}
\setcounter{Corollary}{0} \setcounter{Lemma}{0}
\setcounter{Definition}{0} \setcounter{Remark}{0}
\setcounter{Conjecture}{0}  \setcounter{Example}{0} \hspace{4mm}

We now define supervertices for a $k$-uniform hypergraph.
%-----------------------------------------
\begin{Definition}\label{def-super}
Let $G=(V,E)$ be a $k$-uniform hypergraph.   Let $i \in V$.  The
vertex set
$$U = \{ j \in V : E_j = E_i \}$$
is called a supervertex of $G$.   Clearly, any vertex in the same
supervertex has the same degree.   We call this degree the degree of
that supervertex.   In particular, if a supervertex contains a core
vertex, then all vertices in that supervertex are core vertices.  We
call such a supervertex a core supervertex.  Otherwise, we call it
an intersection supervertex.
\end{Definition}

For example, in Figure 1, there are four intersection supervertices
$\{ 1, 2, 3 \}, \{ 4, 5, 6 \}, \{ 7, 8, 9 \}$, $\{ 10, 11, 12 \}$.
For a loose cycle or a generalized loose $s$-cycle with $1 \le s <
{k \over 2}$, there are $m$ core supervertices and $m$ intersection
supervertices, where $m$ is the number of edges in that $s$-cycle.
Each core supervertex has cardinality $k-2s$. Each intersection
supervertex has degree $2$ and cardinality $s$.

We now have the following theorem.

\begin{Theorem} \label{t5}
Suppose that $G= (V, E)$ is a $k$-uniform hypergraph with $k \ge 3$.
Let $U$ be a supervertex of $G$, with degree $d$ and cardinality
$|U| \ge 2$. Suppose that $\lambda$ is a Laplacian H-eigenvalue of
$G$, $\lambda \not = d$.   Let $\bf x$ be a Laplacian H-eigenvector
of $G$, corresponding to $\lambda$.   Suppose $i, j \in U$.  Then
$|x_i| = |x_j|$.   If $k$ is odd, then $x_i = x_j$.
\end{Theorem}
\noindent {\bf Proof.} By (\ref{e-1}), we have
$$\lambda x_i^{k-1} = dx_i^{k-1} - x_j\sum_{e \in U}\Pi_{s \in e \setminus \{ i, j \}}x_s,$$
and
$$\lambda x_j^{k-1} = dx_j^{k-1} - x_i\sum_{e \in U}\Pi_{s \in e \setminus \{ i, j \}}x_s.$$
Thus,
$$(\lambda -d)x_i^k = (\lambda -d)x_j^k.$$
As $\lambda \not = d$, we have $x_i^k = x_j^k$.   The conclusions
follow from this equality.
 \ep

 Note that Lemma 3.1 of \cite{hqs13} is a special case of this
 theorem.

%%%%%%%%%%%%%%%%%%%%%%%%%%%%%%%%%%%%%%%%%%%%%%%%%%%%%%%%%%%%%%%%%%%%%%%%%%%%%%%
\section{Odd-Uniform Generalized Loose $s$-Cycles}\label{sec-ph2}
\setcounter{Theorem}{0} \setcounter{Proposition}{0}
\setcounter{Corollary}{0} \setcounter{Lemma}{0}
\setcounter{Definition}{0} \setcounter{Remark}{0}
\setcounter{Conjecture}{0}  \setcounter{Example}{0} \hspace{4mm}

In this section, assuming that $k$ is odd, using Theorem \ref{t5},
we show that the largest Laplacian H-eigenvalue of an odd-uniform
generalized loose $s$-cycle is equal to $2$, the maximum degree of
that $s$-cycle.   This result extends the result on odd-uniform
loose cycles in Subsection 4.1 of \cite{hqs13}. Since $k$ is odd,
the case that $k=2s$ is not included. Thus, we have $1 \le s < {k
\over 2}$. The $s$-cycle has always core vertices.

\begin{Proposition} \label{p5}
Suppose that $G= (V, E)$ is an $s$-cycle with $1 \le s < {k \over
2}$ and $k$ is odd. Then $\lambda(\mathcal L) = \Delta = 2$.  If $s$
is even, then the only Laplacian H-eigenvalue $\lambda$ of $G$,
satisfying $\lambda > 1$, is $2$.
\end{Proposition}
\noindent {\bf Proof.}  Suppose that $G$ has $m$ edges.  Then $G$
has $m$ core supervertices $V_i, i \in [m]$ and $m$ intersection
supervertices $U_i, i \in [m]$, displayed as $U_1, V_1, U_2, V_2,
\cdots, U_m, V_m, U_{m+1} \equiv U_1$, such that the edges of $G$
are $e_i = U_i \cup V_i \cup U_{i+1}, i \in [m]$.   For $i \in [m],
|U_i| = s, |V_i|=k-2s$.   Furthermore, assume that the vertices of
$G$ are $j \in [n]$, where $n=m(k-s)$.   Then $U_i = \{ (i-1)(k-s)+j
: j \in [s] \}$, $V_i = \{ (i-1)(k-s)+s+j : j \in [k-2s]$ \}, for $i
\in [m]$.

Suppose that $\lambda > 1$, $\lambda \not = 2$ is a Laplacian
H-eigenvalue of $G$. Let $\bf z$ be an Laplacian H-eigenvector
corresponding to $\lambda$. By Theorem \ref{t5}, we may assume that
for $i \in m$, $y_i = z_{j}$ for $j \in U_i$, $x_i =
z_{j}$ for $j \in V_i$.  Let $y_{m+1} \equiv y_1$,
$x_{0} \equiv x_m$, $y_{0} \equiv y_m$.

By (\ref{e-1}), for $i \in [m]$, we have
\begin{equation} \label{e7.1}
\lambda x_i^{k-1} = x_i^{k-1} - x_i^{k-2s-1}y_i^sy_{i+1}^s
\end{equation}
and
\begin{equation} \label{e7.2}
\lambda y_i^{k-1} = 2y_i^{k-1} - x_i^{k-2s}y_i^{s-1}y_{i+1}^s -
x_{i-1}^{k-2s}y_i^{s-1}y_{i-1}^s.
\end{equation}
If $s$ is even, by (\ref{e7.1}), we must have $x_i = 0$ for all $i \in
[m]$, otherwise we would have $\lambda x_i^{k-1}>x_i^{k-1}> x_i^{k-1} - x_i^{k-2s-1}y_i^sy_{i+1}^s$ for some $i$, contradicting (\ref{e7.1}).  This implies that $y_i \not = 0$ for at least one $i$.   By
(\ref{e7.2}), this implies that $\lambda = 2$, a contradiction. This
proves that when $s$ is even, $\lambda >1$ implies $\lambda = 2$.
The conclusion for the case that $s$ is even is proved.

From now we assume that $s$ is odd and $\lambda > 2$. Then by
(\ref{e7.1}), we see that $x_i \not = 0$ implies that $y_iy_{i+1} <
0$.

(i). First assume that $x_i \not = 0$ for all $i \in [m]$.   By
(\ref{e7.1}), we have $y_iy_{i+1} < 0$ for $i \in [m]$.  Thus, $m$
must be even, otherwise we get a contradiction by the rule of alternating signs of $y_1,\cdots,y_m$.
Assume that $m$ is even.  By (\ref{e7.2}), we have
\begin{equation} \label{e7.3}
(\lambda -2)y_i^k = - x_i^{k-2s}y_i^sy_{i+1}^s -
x_{i-1}^{k-2s}y_{i-1}^sy_i^s.
\end{equation}
Then we have
$$y_1^k - y_2^k + y_3^k - \cdots +y_{m-1}^k - y_m^k = 0.$$
Since $y_iy_{i+1} < 0$ for $i \in [m]$, we get a contradiction, as
$y_1^k - y_2^k + y_3^k - \cdots +y_{m-1}^k - y_m^k$ should have the
same sign as $y_1$, which is nonzero.  The conclusion follows now.

(ii). We now assume $x_i = 0$ for all $i \in [m]$.  This implies
that $y_i \not = 0$ for at least one $i$.   By (\ref{e7.2}), this
implies that $\lambda = 2$, a contradiction.

(iii). Finally, we assume that $x_i = 0$ for some $i \in [m]$ and
$x_i \not = 0$ for other $i \in [m]$.   Without loss of generality,
we may assume that $x_1 > 0$ and $x_m = 0$.  By (\ref{e7.1}), we
have $y_1y_2 < 0$.   By taking $i=1$ in (\ref{e7.3}), we have
$$(\lambda -2)y_1^k = -x_1^{k-2s}y_1^sy_2^s.$$
This implies that $y_1 > 0$.

Now we use induction to show that $x_iy_i>0$ for $i=1,\cdots, m$. The case $i=1$ follows from $x_1>0$ and $y_1>0$. We assume $i\ge 2$ and $x_{i-1}y_{i-1}>0$. Then $y_i\ne 0$ since $x_{i-1}\ne 0$ implying $y_{i-1}y_i<0$. From (\ref{e7.3}) we have (since $k$ and $s$ are both odd):
\begin{equation} \label{e7.4}
0 < (\lambda-2)y_i^{k-s} + x_{i-1}^{k-2s}y_{i-1}^{s}= -x_i^{k-2s}y_{i+1}^s,
\end{equation}
which implies $x_i\ne 0$ and $x_iy_{i+1}<0$. But $x_i\ne 0$ also implies $y_iy_{i+1}<0$, so we obtain $x_iy_i>0$ and thus complete the inductive proof.
Taking $i=m$ in $x_iy_i>0$, we obtain $x_m\ne 0$, a contradiction.

Thus, when $s$ is odd, we cannot have $\lambda > 2$.   This implies
that $\lambda(\mathcal L) = \Delta = 2$. \ep

%%%%%%%%%%%%%%%%%%%%%%%%%%%%%%%%%%%%%%%%%%%%%%%%%%%%%%%%%%%%%%%%%%%%%%%%%%%%%%%
\section{Odd-Uniform Tight $s$-Cycles}\label{sec-ph3}
\setcounter{Theorem}{0} \setcounter{Proposition}{0}
\setcounter{Corollary}{0} \setcounter{Lemma}{0}
\setcounter{Definition}{0} \setcounter{Remark}{0}
\setcounter{Conjecture}{0}  \setcounter{Example}{0} \hspace{4mm}

In this section, we assume that $k$ is odd and $s=k-1$.  Then we
have tight $s$-cycles.   We will see that the results on the largest
Laplacian H-eigenvalues here are very different from those in the
last section.

\begin{Proposition} \label{p6}
Suppose that $G= (V, E)$ is a tight $s$-cycle with $s = k-1$ and $k=
4l+3$ for a nonnegative integer $l$. Then $\Delta = k$.  When $n$,
the number of vertices, is even, we have $\lambda(\mathcal L) \ge
\Delta +1 = k+1$.
\end{Proposition}
\noindent {\bf Proof.} Let ${\bf x} \in \Re^n$ be defined by
$x_{2i-1}= 1$ and $x_{2i} = -1$ for $i \in \left[{n \over
2}\right]$.   Also, we assume that $x_{n+i} \equiv x_i$ for any $i$.
From this we see that the sum of any $k-1$ consecutive components of
$\bf x$ is zero, and the product of any $k-1=4l+2$ consecutive
components of $\bf x$ is $(-1)^{2l+1} = -1$.  Thus we have
$$\prod_{t=j}^{j+k-1}x_t=x_j\prod_{t=j+1}^{j+k-1}x_t=-x_j, \qquad  \mbox {and}  \qquad  \sum_{j=i-k+1}^ix_j=\sum_{j=i-k+1}^{i-1}x_j+x_i=x_i=x_i^k $$
Now multiplying the both sides of (\ref{e-1}) by $x_i$, we obtain
$$\lambda x_i^{k} = kx_i^{k} -
\sum_{j=i-k+1}^i\prod_{t=j}^{j+k-1}x_t= kx_i^{k} +
\sum_{j=i-k+1}^ix_j= kx_i^{k} +x_i^k $$

This shows that $\lambda = k+1$ and $\bf x$ satisfy this system,
i.e., $\lambda = k+1$ is an H-eigenvalue of $\cal L$.  As $\Delta =
k$, the conclusion follows. \ep

This is the second example that $\lambda(\mathcal L) > \Delta$ when
$k$ is odd.   The first example for this is the $3$-uniform complete
hypergraph, given in \cite{hqx13}.  By using supervertices, we may
generalize this result to $k$-uniform regular $s$-cycles, with
$k=q(k-s)$, $q=4l+3$ for some $l$, where $k-s$ and $m$ are even.

We do not know what kind of result can be established for $k=4l+1$.
But for $k=3$, we can get the exact value of $\lambda(\mathcal L)$
when $n=4$, and an upper bound of $\lambda(\mathcal L)$ for all $n$.

\begin{Proposition} \label{p7}
Suppose that $G= (V, E)$ is a tight $s$-cycle with $s = 2$ and
$k=3$. Then $\Delta = 3$.  When $n=4$, we have $\lambda(\mathcal L)=
4$. When $n \ge 5$, we have $\lambda(\mathcal L) \le \Delta +1.5 =
4.5$.
\end{Proposition}
\noindent {\bf Proof.} When $k=3, s=2$ and $n=4$, (\ref{e-1}) has
the form:
\begin{equation} \label{8.1} \left\{\begin{array} {rcl}
(\lambda - 3)x_1^2 & = & - x_2x_3 - x_2x_4 - x_3x_4,\\
(\lambda - 3)x_2^2 & = & - x_1x_3 - x_1x_4 - x_3x_4,\\
(\lambda - 3)x_3^2 & = & - x_1x_2 - x_1x_4 - x_2x_4,\\
(\lambda - 3)x_4^2 & = & - x_1x_2 - x_1x_3 -
x_2x_3.\end{array}\right.
\end{equation}
Summing up these four equations, we have
$$(\lambda - 3)\sum_{i=1}^4 x_i^2 = -2\sum_{1 \le i < j \le 4}x_ix_j
= \sum_{i=1}^4 x_i^2 - (x_1+x_2+x_3+x_4)^2 \le \sum_{i=1}^4 x_i^2.$$
Since $\sum_{i=1}^4 x_i^2 > 0$, we have $\lambda \le 4$.   Combining
with the conclusion of Proposition \ref{p6}, we see that $\lambda(\mathcal L) = 4 =
\Delta +1$.

When $k=3, s=2$, and $n \ge 5$, (\ref{e-1}) has the form:
$$(\lambda - 3)x_i^2 = - x_{i-2}x_{i-1} - x_{i-1}x_{i+1} - x_{i+1}x_{i+2},$$
for $i \in [n]$.  Summing up it for $i$ from $1$ to $n$, we have
$$(\lambda - 3)\sum_{i=1}^n x_i^2 = -\sum_{i=1}^n \left(x_{i-1}x_i + x_{i-1}x_{i+1}
+ x_ix_{i+1}\right) \le {1 \over 2}\sum_{i=1}^n
\left(x_{i-1}^2+x_i^2+x_{i+1}^2\right) = {3 \over 2}\sum_{i=1}^n
x_i^2.$$ Since $\sum_{i=1}^n x_i^2 > 0$, we have $\lambda \le 4.5$.
Thus, $\lambda(\mathcal L) \le 4.5$ in this case.
 \ep

Using Matlab to solve (\ref{8.1}), we find that when $k=3, s=2$ and
$n=4$, the $2$-cycle has only three distinct H-eigenvalues: $4$, $3$
and $0$.

The second conclusion of Proposition \ref{p7} is not sharp in the
proof.  Actually, our Matlab computation shows that when $k=3, s=2$
and $n=5$, the $2$-cycle has only three distinct H-eigenvalues: $3$,
$2.3966$ and $0$; when $k=3, s=2$ and $n=6$, the $2$-cycle has only
four distinct H-eigenvalues: $4$, $3$, $1.7401$ and $0$. Thus, we
have the following conjecture.

\begin{Conjecture} \label{coj1}
Suppose that $k=3$ and $s=2$.   Then $\Delta = 3$.  When $n$ is
even, we have $\lambda(\mathcal L) = 4$.   When $n$ is odd, we have
$\lambda(\mathcal L) = 3$.
\end{Conjecture}

%%%%%%%%%%%%%%%%%%%%%%%%%%%%%%%%%%%%%%%%%%%%%%%%%%%%%%%%%%%%%%%%%%%%%%%%%%%%%%%
\section{Final Remarks}\label{sec-ph4}
\setcounter{Theorem}{0} \setcounter{Proposition}{0}
\setcounter{Corollary}{0} \setcounter{Lemma}{0}
\setcounter{Definition}{0} \setcounter{Remark}{0}
\setcounter{Conjecture}{0}  \setcounter{Example}{0} \hspace{4mm}

In this paper, we showed that the largest signless Laplacian
H-eigenvalue of a connected $k$-uniform hypergraph $G$, where $k \ge
3$, reaches its upper bound $2\Delta$, where $\Delta$ is the largest
degree of $G$, if and only if $G$ is regular, and that the largest
Laplacian H-eigenvalue of $G$, reaches the same upper bound, if and
only if $G$ is regular and odd-bipartite.    We proved that an
even-uniform $s$-path and an even-uniform non-regular $s$-cycle are
always odd-bipartite.  Theorem 4.1 characterized odd-bipartite
regular $s$-cycles.   We identified the largest signless Laplacian
H-eigenvalue of an $s$-cycle.   When the $s$-cycle is odd-bipartite,
this gives the largest Laplacian H-eigenvalue of that $s$-cycle. We
then introduced supervertices and showed that the largest Laplacian
H-eigenvalue of an odd-uniform generalized loose $s$-cycle is $2$,
the maximum degree of that $s$-cycle.   We also showed that the
largest Laplacian H-eigenvalue of a $k$-uniform tight $s$-cycle is
not less than the maximum degree of that $s$-cycle, plus one, if the
number of vertices is even and $k=4l+3$.  It will be a further
research topic to prove or to disprove Conjecture 8.1, and to
identify  the largest Laplacian H-eigenvalue of an $s$-path or a
general non-odd-bipartite $s$-cycle, for $s \ge 2$. It will be
interesting to see if one may use the tensor eigenvalue theory to
study other research topics related with $s$-paths and $s$-cycles.

%{\bf Acknowledgement.}  We are thankful to the referee for his or
%her comments, which helped us to improve our paper.
%%%%%%%%%%%%%%%%%%%%%%%%%%%%%%%%%%%%%%%%%%%%%%%%%%%%%%%%%%%%%%%%%%%%%%%%%%%%%%%%%%%%%%%%%%%%%%
\bibliographystyle{model6-names}

\begin{thebibliography}{00}
%\bibitem{b73}
%Berge, C., 1973. Hypergraphs. Combinatorics of finite sets, third edition. North-Holland, Amsterdam.

\bibitem{bh11}
Brouwer, A.E., Haemers, W.H., 2011. Spectra of Graphs. Springer. New
York.

\bibitem{cqz13}
Chang, K.C., Qi, L., Zhang, T., 2013. A survey on the spectral
theory of nonnegative tensors. Numerical Linear Algebra with
Applications, in press.


%\bibitem{c97}
%Chung, F.R.K., 1997. Spectral Graph Theory. Am. Math. Soc.

%\bibitem{cd12}
%Cooper, J., Dutle, A., 2012. Spectra of uniform hypergraphs. Linear Algebra Appl., 436, 3268--3292.

%\bibitem{er60}
%Erd\"os, P., Rado, R., 1960.  Intersection theorems for systems of
%sets. Journal of the London Mathematical Society, Second Series 35,
%85--90.

\bibitem{f10}
Frieze, A., 2010. Loose Hamilton cycles in random 3-uniform
hypergraphs, Electron. J. Combin. 17, \#28.

\bibitem{fkl12}
Frieze, A., Krivelevich, M., Loh, P.-S., 2012.  Packing tight
Hamilton cycles in 3-uniform hypergraphs. Random Structures
Algorithms, 40, 269--300.

%\bibitem{gm94}
%Grone, R., Merris, R., 1994. The Laplacian spectrum of a graph. II.
%SIAM J. Discrete Math., 7, 221--229.

\bibitem{gss08}
Gy\'arf\'as, A., S\'ark\"ozy, G., Szemer\'edi, E., 2008. The Ramsey
number of diamond-matchings and loose cycles in hypergraphs.
Electronic Journal of Combinatorics, 15, \#R126.

\bibitem{hlprrss06}
Haxell, P., Luczak, T., Peng, Y., R\"odl, V., Ruci\'nski, A.,
Simonovits, M., Skokan, J., 2006. The Ramsey number for hypergraph
cycles I.  J. Combin. Theory Ser. A., 113, 67--83.

%\bibitem{hhlq12}
%Hu, S., Huang, Z.H., Ling, C., Qi, L., 2013. On determinants and eigenvalue theory of tensors.
%J. Symb. Comput., 50, 508--531.

%\bibitem{hq12a}
%Hu, S., Qi, L., 2012. Algebraic connectivity of an even uniform
%hypergraph. J. Comb. Optim., 24, 564--579.

%\bibitem{hq13}
%Hu, S., Qi, L., 2013. The Laplacian of a uniform hypergraph. J. Comb. Optim., in press.

\bibitem{hq13}
Hu, S., Qi, L., 2013. The eigenvectors of the zero Laplacian and signless Laplacian eigenvalues of a uniform hypergraph. arXiv:1303.4048.

\bibitem{hqs13}
Hu, S., Qi, L., Shao, J., 2013. Cored hypergraphs, power hypergraphs
and their Laplacian eigenvalues. Linear Algebra and Its
Applications, in press, DOI: 10.1016/j.laa.2013.08.028.

\bibitem{hqx13}
Hu, S., Qi, L., Xie, J., 2013. The largest Laplacian and signless Laplacian H-eigenvalues of a uniform hypergraph. arXiv:1304.1315.

\bibitem{kkmo11}
Keevash, P., K\"uhn, D., Mycroft, R., Osthus, D., 2011. Loose
Hamilton cycles in hypergraphs. Dis. Math., 311, 544--559.

\bibitem{kmo10}
K\"uhn, D., Mycroft, R., Osthus, D., 2010. Hamilton $l$-cycles in
uniform hypergraphs. J. Combin. Theory Ser. A., 117, 916--927.


\bibitem{ko06}
K\"uhn, D., Osthus, D., 2006. Loose Hamilton cycles in 3-uniform
hypergraphs of high minimum degree. J. Combin. Theory Ser. B., 96,
767--821.

\bibitem{lqy12}
Li, G., Qi, L., Yu, G., 2013. The Z-eigenvalues of a symmetric tensor and its application to spectral hypergraph theory.
Numer. Linear Algebr., in press.

\bibitem{mors13}
Maherani, L., Omidi, G.R., Raeisi, G., Shasiah, J., 2013. The Ramsey
number of loose paths in 3-uniform hypergraphs. Electron. J.
Combin., 20, \#12.

%\bibitem{n13}
%Nikiforov, V., 2013. An analytical theory of extremal hypergraph
%problems. arXiv:1305.1073.

%\bibitem{pz12}
%Pearson, K.J., Zhang, T., 2013. On spectral hypergraph theory of the
%adjacency tensor. Graphs and Combinatorics, in press.

\bibitem{p13}
Peng, X., 2013. The Ramsey number of generalized loose paths in
uniform hypergraphs. arXiv:1305.0294.

\bibitem{q05}
Qi, L., 2005. Eigenvalues of a real supersymmetric tensor. J. Symb.
Comput., 40, 1302--1324.

\bibitem{q13}
Qi, L., 2013. H$^+$-eigenvalues of Laplacian and signless Laplacian
tensors. Communications in Mathematical Sciences, in press.
 arXiv:1303.2186.

%\bibitem{q12b}
%Qi, L., 2013. Symmetric nonnegative tensors and copositive tensors.
%Linear Algebra Appl., 439, 228-238.

%\bibitem{r09}
%Rota Bul\`{o}, S., 2009. A game-theoretic framework for similarity-based
%data clustering. PhD thesis, Universit\`{a} Ca' Foscari di Venezia,
%Italy.

%\bibitem{rp09}
%Rota Bul\`{o}, S., Pelillo, M., 2009. A generalization of the
%Motzkin-Straus theorem to hypergraphs. Optim. Lett., 3, 187--295.

%\bibitem{sqh13}
%Shao, J., Qi, L., Hu, S., 2013. Some new trace formulas of tensors
%with applications in spectral hypergraph theory. arXiv:1307.5690.

%\bibitem{s03}
%Stevanovi\'c, 2003. Bounding the largest eigenvalue of trees in
%terms of the largest vertex degree, Linear Algebra Appl., 360,
%35-42.

%\bibitem{xc12c}
%Xie, J., Chang, A., 2013. H-Eigenvalues of the signless Laplacian
%tensor for an even uniform hypergraph. Front. Math. China, 8,
%107-127.

%\bibitem{xc12a}
%Xie, J., Chang, A., 2012. On the Z-eigenvalues of the signless Laplacian tensor for an even uniform hypergraph.
%Manuscript, Fuzhou University.

%\bibitem{yy11}
%Yang, Q., Yang, Y., 2011. Further results for Perron-Frobenius
%theorem for nonnegative tensors II. SIAM J. Matrix Anal. Appl., 32,
%1236-1250.

\bibitem{z07}
Zhang, X.D., 2007, The Laplacian eigenvalues of graphs: a survey.
Gerald D. Ling (Ed.), Linear Algerbra Research Advances, Nova
Science Publishers Inc. Chapter 6, pp. 201--228.
\end{thebibliography}

\end{document}